\documentclass[a4paper,10pt]{amsart}
\usepackage[english]{babel}
\usepackage[utf8]{inputenc}
\usepackage[T1]{fontenc}
\usepackage{csquotes} 
\usepackage[backend=bibtex, style=numeric,
    useprefix=false,%
	giveninits=true,%
	doi=false,%
	url=true,%
	isbn=false,%
	maxbibnames=99%
	]{biblatex}
\bibliography{./BIB-JetSpacesOnCarnotGroups}
\usepackage[normalem]{ulem}

\usepackage{amssymb}
\usepackage{mathrsfs}
\usepackage[all]{xy}	
\usepackage{hyperref}
\usepackage[usenames,dvipsnames]{xcolor}
\hypersetup{colorlinks,%
citecolor=Black,%
filecolor=Black,%
linkcolor=Black,%
urlcolor=Black}
\usepackage{enumitem}	
\usepackage{array}		
\allowdisplaybreaks
\usepackage{graphicx}
\usepackage{pdflscape}

\newcommand{\dbd}[2]{\frac{\partial#1}{\partial #2}}

\newcommand{\scr}[1]{\mathscr{#1}}
\newcommand{\subh}[1]{ {#1}_{\scriptscriptstyle \mathscr{H} }}

\newcommand{\frk}[1]{\mathfrak{#1}}

\newcommand{\N}{\mathbb{N}}	
\newcommand{\R}{\mathbb{R}}	
\newcommand{\He}{\mathbb{H}}

\newcommand{\Lin}{\mathtt{Lin}}

\newcommand{\jet}{\mathfrak{j}}
\newcommand{\Jet}{\mathtt{J}}
\newcommand{\HorDer}{\mathtt{HD}}

\newcommand{\Sp}{\operatorname{Sp}}

\newcommand{\rcontr}{\raisebox{\depth}{\scalebox{1}[-1]{$\lrcorner$}}} 

\theoremstyle{plain}
\newtheorem{proposition}{Proposition}[section]
\newtheorem{theorem}[proposition]{Theorem}

\theoremstyle{definition}

\theoremstyle{remark}

\title{Moduli of Sub-Laplacians on the second Heisenberg group }
\author[Nicolussi~Golo]{Sebastiano Nicolussi Golo}
\address[Nicolussi~Golo]{currently unafiliated, Italy}	
\email{sebastiano2.72@gmail.com}

\author[Warhurst]{Benjamin Warhurst}
\address[Warhurst]{Institute of Mathematics, University of Warsaw, ul. Banacha 2, 02-097 Warsaw, Poland}
\email{b.warhurst@mimuw.edu.pl}
\thanks{
 S.~N.~G.~has been supported by the Academy of Finland (%
grant 328846, ``Singular integrals, harmonic functions, and boundary regularity in Heisenberg groups'',
grant 322898 ``Sub-Riemannian Geometry via  Metric-geometry and Lie-group Theory'',
grant 314172 ``Quantitative rectifiability in Euclidean and non-Euclidean spaces''),
and by the University of Padova STARS Project ``Sub-Riemannian Geometry and Geometric Measure Theory Issues: Old and New''.\\
$\text{ }$ $\text{ }$  S.~N.~G.~and B.~W.~are grateful for the support provided by the grant of the National Science Center, Poland (NCN), UMO-2017/25/B/ST1/01955.\\
$\text{ }$ $\text{ }$   We also would like to thank Prof Ian Anderson, Utah State University, for his insight and tips with MAPLE computation.}

\keywords{Jet spaces; stratified Lie groups; Carnot groups; Laplacian; equivqlence; }

\subjclass[2010]{%
58A20; 
22E25; 
35R03; 
}

\date{\today} 

\begin{document}
\begin{abstract}
We solve the contact equivalence problem for generalised sub-Laplacians on $\He^2$ and show that the family of sub-Laplacians on $\He^2$ modulo contact equivalence, is parameterised by $\R^+$.  
\end{abstract}
\maketitle

\setcounter{tocdepth}{2}
\phantomsection
\addcontentsline{toc}{section}{Contents}
\tableofcontents
\section{Introduction}

Each positive definite $A \in GL(\R^n)$ determines the second order operator $$\Delta_A=\nabla \cdot A \nabla = \sum_{ij} a_{ij}  \frac{\partial^2}{\partial x_i \partial x_i }.$$ Since $A$ is positive definite, it has a symmetric square root $C$, which is to say $C$ is symmetric and $C^2 = A^{-1}$. If $f(x)=Cx$, then $$\Delta_A (u \circ f) = (\Delta u) \circ f. $$

Let $G$ be a connected, simply connected, stratified nilpotent Lie group, with Lie algebra $\frak{g} = \mathbb{V}_1 \oplus \dots \oplus \mathbb{V}_s$, where $\mathbb{V}_{i+1} = [\mathbb{V}_1, \mathbb{V}_i]$ and $[\frak{g}, \mathbb{V}_s]=0$. Let $d_i = {\rm dim}\, \mathbb{V}_i$ and $N={\rm dim} \, \frak{g}$. If $e_i$ denotes an element of the standard basis of $\R^N$, then we choose $ \mathbb{V}_1={\rm span}\{ e_1,\dots, e_{d_1} \}$,   $ \mathbb{V}_2={\rm span}\{ e_{d_1+1},\dots, e_{d_2} \}$ and so on.  The left invariant vector fields are framed by the set $\scr{F}=\{\tilde X_1, \dots, \tilde X_N\}$ where $$\tilde X_i u(p) = \frac{d}{dt} u( p\exp(t e_i)) |_{t=0}$$ and the horizontal bundle, denoted by $\scr{H}$, is the subbundle of $TG$ framed by  $\{\tilde X_1 \dots, \tilde X_{d_1}\}$. The quadratic form on $\scr{H}$ which makes $\{\tilde X_1, \dots, \tilde X_{d_1} \}$ orthonormal is left invariant and denoted by $g_0$. If $\tilde V$, $\tilde W$ are $\scr{H}$ valued left invariant vector fields, with coefficients $v_i$ and $w_i$ relative to $\scr{F}$, then we have that $$g_0( \tilde V, \tilde W) =( v_1, \dots,  v_{d_1}) \cdot ( w_1, \dots,  w_{d_1}).$$  

A local diffeomorphism $F: U \subset G \to G$, where $U \subset G$ is open, which preserves the horizontal bundle in the sense that $dF_p(\scr{H}_p) =\scr{H}_{F(p)}$, will be referred to as a contact mapping, despite the fact that the distributions are not necessarily associated with contact structures.  

The horizontal gradient of a real valued function $u$ on $G$ is define to be the $\scr{H}$-valued vector field $$ \subh{\nabla}u= \sum_{i=1}^{d_1} (\tilde X_iu) \tilde X_i. $$ In coordinates we have $\subh{\nabla}=(\tilde X_1, \dots, \tilde X_{d_1})$ and the canonical sub-Laplacian is defined to be $$\mathcal{L}u = \subh{\nabla} \cdot \subh{\nabla}u =  \sum_{i=1}^{d_1}  \tilde X_i^2u.$$
Given a positive-definite symmetric matrix $M=(m_{ij}) \in GL(\R^{d_1})$ we define a generalised sub-Laplacian  $\mathcal{L}_M$ by setting
$$ \mathcal{L}_Mu  = \subh{\nabla} \cdot M  \subh{\nabla} u = \sum_{ij=1}^{d_1} m_{ij}  \tilde X_i   \tilde X_j u.$$

One can then ask if there exists a local contact diffeomorphism $f:G \to G$ such that 
\begin{align}
\mathcal{L}_M(u \circ f) = \mathcal{L}u \circ f. \label{Linvariance}
\end{align}
For free groups, all sub-Laplacians are equivalent, as shown in \cite{BLUstraGrp} (Theorem 16.1.2 page 625) where matters are simplified by the fact that any element of  $GL(\R^{d_1})$ extends to a Lie algebra automorphism. The first Heisenberg group $\He^1$ is free, and so $\He^2$ is perhaps the simplest case not covered by the results in \cite{BLUstraGrp}. In this study, we deal with the equivalence problem for sub-Laplacians by applying the theory of jet spaces over stratified groups, see \cite{SebaBenRevista}. 
 
 \section{Horizontal jets over stratified groups}
 
 The horizontal Taylor expansion of a smooth real valued function $u$ at a point $p \in G$, is based on the usual Taylor expansion around zero of the function $h(t)= u(p\exp(t V))$ when $V \in \mathbb{V}_1$. If $V=\sum_{i=1}^{d_1} v_i e_i$ and $$ \tilde V^k u(p) = \frac{d^k}{dt^k} u( p\exp(t V)) |_{t=0}$$then the $m$-th order Taylor polynomial of $h$ around zero, denoted  $T^m_0 h$, has the form 
 \begin{align*}
 T^m_0 h(t) &= \sum_{k=0}^m \tilde V^ku (p) \frac{t^k}{k!} \\
 &= \sum_{k=0}^{m} \, \sum_{i_1, \dots i_{k}=1}^{d_1} v_{i_1}\dots v_{i_k} \tilde X_{i_1} \dots \tilde X_{i_k} u(p)  \frac{t^k}{k!}.
 \end{align*}
 Setting $t=1$ determines the horizontal Taylor polynomial $T^m_pu$ of order $m$, defined by
 \begin{align}
 T^m_pu(V)
 	&= \sum_{k=0}^m  \frac{1}{k!} \tilde V^ku (p) \nonumber \\
 	&= \sum_{k=0}^{m} \frac{1}{k!} \sum_{ i_1, \dots, i_{k}=1}^{d_1} \, \,   \, \tilde X_{i_1} \dots \tilde X_{i_k} u(p) \, v_{i_1}\dots v_{i_k}, \label{taylormjet}
 \end{align} which approximates $V \to u( p\exp(V))$, $V \in \mathbb{V}_1$.

The $k$-th term of $T^m_pu$ is a homogeneous polynomial of degree $k$, and can be expressed in "polar" type form, by which we mean  $T^m_pu(V)= \sum_{k=0}^m \frac{1}{k!} A^k_{u,p}(V, \dots, V)$, where for $k>0$, $A^k_{u,p}$ is the $k$-linear map defined on $\mathbb{V}_1^k$ by 
\begin{align} A^k_{u,p} \left({ {V} }_{1},\dots ,{ {V} }_{k}\right)
&=  \tilde V_k \dots  \tilde V_1 u(p) \label{polarform}
 \end{align} and $A^0_{u,p}=u(p)$. 
The horizontal $m$-jet of $u$, denoted $j^m_pu$, is the formal sum
$$ j^m_pu= \sum_{k=0}^m \frac{1}{k!} A^k_{u,p}.$$ The horizontal $m$-jet of a vector valued function $u:G \to \mathbb{W}$ is then {\tiny }given by $(j^m_pu_1, \dots, j^m_pu_N)$ where $u_1, \dots , u_N$ are the coordinate functions relative to a basis of the target vector space $\mathbb{W}$.

For $k \in \mathbb{N}$, we denote by $\HorDer^k(\frk g; \mathbb{W})$ the vector subspace of $\Lin^k(\mathbb{V}_1; \mathbb{W})$ consisting of all $k$-multilinear maps of the form $A^k_{u,p}$.
 For $k=0$, we set $\HorDer^0(\frk g; \mathbb{W})= \mathbb{W}$. For each positive integer $m$, we define   
 \[
 \HorDer^{\leq m}(\frk g; \mathbb{W})=\bigoplus_{k=0}^m \HorDer^k(\frk g; \mathbb{W}).  \quad (\HorDer^{\leq m} \text{ for short})
 \]
 The $m$-th  order horizontal jet space, denoted by $\Jet^m(G, \mathbb{W})$, is the vector bundle $G \times \HorDer^{\leq m} \to G$. 
 
 The Lie group structure on $\Jet^m(G, \mathbb{W})$ is given by an anti-semidirect product which we now briefly review. If $\mathfrak{g}$ and $\mathfrak{h}$ are Lie algebras, then an antimorphism $\psi:\mathfrak{g} \to {\rm Der}( \mathfrak{h})$ is a linear map with the property that $$\psi([x,y])=-[\psi(x),\psi(y)]).$$
 An antimorphism induces a Lie bracket on $\mathfrak{g} \times \mathfrak{h}$ by the formula $$[(x,X), (y,Y)]=([x,y], -[X,Y]+\psi(y)X - \psi(x)Y). $$ We denote this Lie algebra by  $\mathfrak{g}\ltimes_\psi \mathfrak{h}$ (anti-semidirect product). 
 
 Let $G$ and $H$ be the  connected simply connected Lie groups corresponding to
 $\mathfrak{g}$ and $\mathfrak{h}$. The antimorphism $\psi$ induces an antimorphism $\phi : G \to {\rm Aut} (H)$, by which we mean \begin{align*}\phi(ab) &= \phi(b)\phi(a) \\ d\phi|_{\exp(x)} &= e^{\psi(x)} \in {\rm Aut} (\mathfrak{h}) \quad  \forall x \in \mathfrak{g}. \end{align*} 
 
 The set  $G \times H$ and the product given by $$(a,A)(b,B) = (ab , B\phi(b)A )$$ 
 define a Lie group $G\ltimes_\phi H$ with Lie algebra $ \mathfrak{g}\ltimes_\psi \mathfrak{h}$. When $\mathfrak{h}$ is abelian, the algebra and group products simplify to the following:  
 \begin{align}
 \nonumber
 [(x,X),(y,Y)] &= ([x,y], \psi(y)X-\psi(x)Y ) , \\
 \nonumber
 (a,A)(b,B) &= (ab , B+\phi(b)A ) .
 \end{align}
 If $(a,A)\in G\ltimes_\phi H $ and $(x,X)\in \mathfrak{g}\ltimes_\psi \mathfrak{h} $, then  the left-invariant vector field generated by $(x, X)$  
 evaluated at $(a,A)\in G\ltimes_\phi H $ is
 \[
 \widetilde{(x,X)} (a,A) = ( \tilde x(a) , \psi(x)A + X ) 
 \]where $\tilde x(a)$ is the left-invariant vector field generated by $x$  evaluated at $a$.

 If $k \geq 2$ and  $A\in \HorDer^k(\frk g ;\mathbb{W})$, then for $V\in \mathbb{V}_1$, the right contraction of $A$ by $V$, denoted $V\rcontr A$, is the element in $\HorDer^{k-1}(\frk g;\mathbb{W})$ given by
 $$
 V\rcontr A (V_1,\dots,V_{k-1}) =  \, A(V_1,\dots,V_{k-1},V).
 $$
 Moreover, if $k = 1$ then $V\rcontr A=A(V)$, and if $k=0$ then $V\rcontr A=0$. 
 
 Given $V\in \mathbb{V}_1$, we extend the definition of the right contraction by $V$ to a map 
 $V\rcontr:\HorDer^{\le m} \to \HorDer^{\le  m }$,
 where for each $A\in\HorDer^{\le m}$, we set 
 $(V\rcontr A)^k = V\rcontr A^{k+1}$ for $k<m$ and $(V\rcontr A)^m=0$. The map $\mathbb{V}_1\to {\rm End}(\HorDer^{\le m})$, $V\mapsto V\rcontr$, extends uniquely (by formula \eqref{polarform}) to a Lie algebra anti-morphism $\psi: \frk g\to {\rm End}(\HorDer^{\le m})$, i.e., 
\begin{align}
[V,W]\rcontr = - [V\rcontr,W\rcontr]. \label{rcontanti}
\end{align}

Since $\HorDer^{\le m}$ is an abelian Lie algebra, ${\rm Der}(\HorDer^{\le m})={\rm End}(\HorDer^{\le m})$. Hence, $\psi(V)= V\rcontr$ is a Lie algebra antimorphism from $\frk g$ to ${\rm Der}(\HorDer^{\le m})$. We obtain a Lie algebra 
\[
\jet^m(\frk g;\mathbb{W}) := \frk g \ltimes_\psi \HorDer^{\le m} ,
\]
where
\begin{align}
[(V,A),(W,B)] = ([V,W],W\rcontr A - V\rcontr B). \label{jLBrkt}
\end{align}
The Lie algebra	$\jet^m(\frk g;\mathbb{W})$ is stratified of step $\max\{s,m+1\}$, where $s$ is the step of $\frk g$.
The $k$-th layer of $\jet^m(\frk g;\mathbb{W})$ is
\[
\jet^m(\frk g;\mathbb{W})_k := \mathbb{V}_k\times \HorDer^{m+1-k}(\frk g;\mathbb{W})
\]
where $\mathbb{V}_k=\{0\}$ if $k>s$ and $\HorDer^{m+1-k}(\frk g;\mathbb{W})=\{0\}$ if $k> m+1$.

The simply connected Lie group corresponding to the Lie algebra $\jet^m(\frk g;\mathbb{W})$ is the semidirect product $G \ltimes \HorDer^{\le m}(\frk g;\mathbb{W})$
with group law
\[
(\exp(V),A)(\exp(W),B) = ( \exp(V)\exp(W) , B+e^{(W\rcontr)}A ). \]
Hence $\Jet^m(G, \mathbb{W})$ together with the afore mentioned product is a stratified Lie group withe stratified Lie algebra $\jet^m(\frk g;\mathbb{W})$. The underlying manifold structure of $\Jet^m(G, \mathbb{W})$ is that of the vector bundle $G \times \HorDer^{\leq m} \to G$. 

 The equivalence problem defined by \eqref{Linvariance} leads to an equivalence by contact maps of two stratified subgroups of the second horizontal jet space over $G=\He^2$. The subgroups are submanifolds $ \Jet \mathcal{L}, \Jet \mathcal{L}_M \subset  \Jet^2(\He^2, \R)$ where  $ \Jet \mathcal{L}$ consists of $2$-jets of $\mathcal{L}$-harmonic functions and  $ \Jet \mathcal{L}_c$ consists of $2$-jets of $\mathcal{L}_c$-harmonic functions, see \eqref{Lc} for the definition of $\mathcal{L}_c$. 
 
 Thus it needs to be established that a contact map of an open set $\Omega \subseteq G$ can be prolonged to a contact map of and open subset of the second horizontal jet space. To this end we show directly that we can prolong from $G$ to $\Jet^1(G, \mathbb{W})$ and then invoke Theorem 5.2 in \cite{SebaBenRevista} to get a prolongation to $\Jet^2(G, \mathbb{W})$.  
 
 When $m=0$ we have $\Jet^0(G, \mathbb{W}) =G \times \mathbb{W}$ where $(a,A)(b,B)=(ab,A+B)$ and $\scr{H}$ is given by the left translation of $\mathbb{V}_1 \times \mathbb{W}$. A contact map $f:\Omega \subseteq G \to G$  lifts to a contact map $F$ on $\Omega \times \mathbb{W}$ by setting $F=(f \circ \pi_G , \pi_{\mathbb{W}})$. 
 Indeed if $(a, A) \in \Jet^0(G, \mathbb{W})=G \times \mathbb{W}$ and $(V, B) \in \jet^0( \frk{g}, \mathbb{W})=\mathbb{V}_1 \times \mathbb{W}$, then $$ \gamma(t)=(a, A)(\exp(tV), tB)=(a \exp(tV), A +tB)$$ represents a horizontal tangent at $(a,A)$. Moreover
 \begin{align*}	
 	F(a,A)^{-1}F \circ \gamma(t) &=(f(a)^{-1}f(a \exp(tV)) , tB)
 \end{align*}
 represents a horizontal tangent at the identity since
 \begin{align*}	
 	\frac{d}{dt} F(a,A)^{-1} F \circ \gamma(t)|_{t=0} &=  \big ( \frac{d}{dt} f(a)^{-1}f(a \exp(tV))|_{t=0} , B \big ) \in \mathbb{V}_1 \times \mathbb{W}
 \end{align*} 
 where $f$ being a contact map implies
 \begin{align}
 N_{f,a}(V):=\frac{d}{dt} f(a)^{-1}f(a \exp(tV))|_{t=0} \in \mathbb{V}_1. \label{Nfa}
 \end{align}

When $m=1$ we have $\Jet^1(G, \mathbb{W}) =G \times \HorDer^{\leq 1}$ where
\begin{align*}
(a,A)(b,B)&=(ab,B+ A+\log(b)\rcontr A^1)
\end{align*} and $\scr{H}$ is given by the left translation of $\mathbb{V}_1 \times \HorDer^{1}$. Moreover, 
 $$(b,B)^{-1}=(b^{-1}, \log(b)\rcontr B^1 - B)$$
 If $\pi_0 :\Jet^1(G, \mathbb{W}) \to \Jet^0(G, \mathbb{W})$ is the natural projection, then a contact map $\hat F$ of $ \tilde U \subseteq \Jet^1(G, \mathbb{W})$ is a prolongation of a contact map $F$ of $\Jet^0(G, \mathbb{W})$ if $\pi_0 \circ \hat F =F \circ \pi_0$, which when we write $\hat F=( \hat F_G, \hat F^0+\hat F^1)$,  implies that $\hat F_G  =F_G \circ \pi_0$ and  $\hat F^0  =F^0 \circ \pi_0$. If we assume $F$ is a prolongation, then $F_G  =f \circ \pi_G \circ \pi_0$ and  $\hat F^0  = \pi_{\mathbb{W}} \circ \pi_0$, in summary
 $$\hat F(a,A)=(f(a), A^0+ \hat F^1(a,A) )$$
 and
 $$\hat F(a,A)^{-1}=(f(a)^{-1}, \log(f(a)) \rcontr \hat F^1(a,A) - A^0 - \hat F^1(a,A) )$$
 If $(a, A) \in \Jet^1(G, \mathbb{W})=G \times  \HorDer^{\leq 1} $ and $(V, B) \in \jet^1( \mathbb{V}_1, \HorDer^{1}) =\mathbb{V}_1 \times \HorDer^{1}$, then $$ \gamma(t)=(a, A)(\exp(tV), tB)=(a \exp(tV), A +t(B+ V \rcontr A) )$$ represents a horizontal tangent at $(a,A)$. It follows that
 \begin{align*}	
 \hat F_G \circ \gamma(t) &=  f(a \exp(tV)) \\
 \hat F^0 \circ \gamma(t) &=  A^0+t \, V \rcontr A^1 \\
 \hat F^1 \circ \gamma(t) &=  \hat F^1 \big ( a \exp(tV), A + t  (B+ V \rcontr A) \big )
\end{align*}
and
\begin{enumerate}	
\item $ \displaystyle \pi_G \big ( \hat F(a, A)^{-1} \hat F \circ \gamma(t) \big ) =  f(a)^{-1} f(a \exp(tV))$
\item $ \displaystyle \pi_{\HorDer^{0}}\big ( \hat F(a, A)^{-1} \hat F \circ \gamma(t) \big ) = -\log \big ( f(a)^{-1} f(a \exp(tV)) \big ) \rcontr \hat F^1 (a,A) +t V \rcontr A^1$ 
\item $ \displaystyle \pi_{\HorDer^{1}}\big ( \hat F(a, A)^{-1} \hat F \circ \gamma(t) \big ) = \hat F^1 \circ \gamma(t)- \hat F^1 (a,A).$
\end{enumerate}
Note that in item (2) we use the fact that if $a=\exp(V)$ and $b=\exp(W)$ then 
$$ \log(a^{-1}b)= \log(a)-\log(b) + H(V,W)$$ where $H$ is an expression involving brackets of order greater than one  which implies
 $$ \log(a) \rcontr \hat F^1(a, A)  -\log(b) \rcontr \hat F^1(a, A)  = -\log \big ( a^{-1} b \big ) \rcontr \hat F^1 (a,A). $$ 
 
If $\hat F$ is to be a contact map, then we require that the expression on the right in item (2) vanishes which implies 
\begin{align}	
 V \rcontr A^1 =N_{f,a}(V) \rcontr \hat F^1 (a,A). \label{detmCond1}
\end{align}

Since $f$ is a contact map, the set $\{N_{f,a}(e_i) : i=1\dots d_1\}$ is a basis for $\mathbb{V}_1$ and $N_{f,a}$ is an invertible linear map of $\mathbb{V}_1$ to itself. It now follows from \eqref{detmCond1}  that $$\hat F^1(a,A)= A^1 \circ N_{f,a}^{-1}.$$ If $$\eta(g)=A^1 \circ N_{f,a}^{-1} \circ \log \circ L_{f(a)}^{-1}(g) + A^0$$
then $\eta(f(a))=A^0$ and $$ \tilde V \eta(f(a)) = A^1 \circ N_{f,a}^{-1}(V),$$ hence if follows that $$\hat F^0(a, A) + \hat F^1(a, A)= A_{\eta, f(a)}^0 + A_{\eta, f(a)}^1$$ is an element of $ \HorDer^{\leq 1}$. The same construction starting with $f^{-1}$, then prolonging to $F^{-1}$, and then prolonging to $\hat F^{-1}$ yields the inverse of $\hat F$ and so $\hat F$ is a diffeomorphism.

 When $m = 2$, we apply Theorem 5.2 in \cite{SebaBenRevista}, see Remark 5.4.

\section{Heisenberg Groups} 

Following \cite{KORANYI19951}, the Heisenberg group $\He^n$ is given by $\R^{2n+1}$, where we denote coordinates by $(x,y,t)$, $x,y \in \R^n$ and $t \in \R$, and the product is given by
$$(x,y,t)(x',y',t')=(x+x', y+y', t+t'-2 x \cdot y'+2 y \cdot x')).$$
If $\{e_i:i=1,\dots, 2n+1\}$ is the standard basis for $\R^{2n+1}$ and $\frak{h}^n$ denotes the Lie algebra where the nontrivial Lie brackets are $[e_i,e_{n+i}]=-4 e_{2n+1}$, $i=1,\dots,n$, then the product above is the Baker--Campbell--Hausdorff model of $\He^2$ built on $\frak{h}^n$.

The left invariant vector fields are framed by
\begin{align*} 
 \tilde X_i &= \dbd{}{x_i}+2y_j \dbd{}{t} \quad i=1,\dots , n\\
 \tilde Y_i &= \dbd{}{x_i}-2x_j \dbd{}{t}  \quad i=1,\dots , n\\
 \tilde T &=\dbd{}{t}
\end{align*}   
for which the only nontrivial brackets are  $$[\tilde X_i,\tilde Y_i] =-4 \tilde T \quad i=1,\dots,n.$$ 

The frame $\{dx_i, dy_i, \theta \}$, where $\theta =dt +2 \sum_i (x_i dy_i -y_i dx_i)$, is the dual of the frame $\{\tilde X_i, \tilde Y_i, \tilde T \}$ and ${\rm ker \, } \theta=\scr{H}$. In particular, $\theta$ is a contact form since $$(d \theta)^n \wedge \theta =2^{2n} dx_1 \wedge \dots \wedge dx_n \wedge dy_1 \wedge \dots \wedge d y_n \wedge dt. $$
The left invariant $2$-form $$ B=\frac{1}{4} d \theta =  \sum_i dx_i \wedge dy_i $$ is symplectic on $\scr{H}$  and represented by the matrix $$\Omega=\begin{pmatrix}
0 & I_n\\
-I_n & 0
\end{pmatrix}$$ relaltive to the chosen coordinates.  Furthermore, Cartan's formula implies that  $$[\tilde V, \tilde W ]= -4 B(\tilde V, \tilde W ) \tilde T$$ for all left invariant horizontal vector fields $\tilde V$ and $\tilde W $. The complex structure $J_0=\Omega^T$ is compatible with $B$ and $g_0$, that is to say $$g_0(\tilde V,\tilde W )=B(\tilde V ,J_0\tilde W )$$ for all left invariant horizontal vector fields $\tilde V$ and $\tilde W$. 
 
A local diffeomorphism $f:\He^n \to \He^n$ preserves $\scr{H}$ if and only $f^*\theta= f^*\theta(\tilde T) \theta$ and is therefore a contact transformation in the strict sense. Left translations are contact transformations.

An automorphism $\alpha$ of $\He^n$ must satisfy $d \alpha (\tilde T)= \lambda \tilde T$ since $\tilde T$ is central, and is therefore a contact transformation. Two fundamental examples are dilation by $\lambda >0$ 
\begin{align*} 
\delta_{\lambda}(x,y,t) &= (\lambda x, \lambda y, \lambda^2 t) 
\end{align*} 
and the reflection 
\begin{align} 
r(x,y,t) &= ( x, - y, -t). \label{reflection}
\end{align}
Consequently, every automorphism of $\He^n$ is a composition $ \alpha \circ \delta_{\lambda} \circ r$ where $\alpha$ is normalised, that is to say, $d \alpha (\tilde T)=\tilde T$. Furthermore, a normalised automorphism $\alpha$ satisfies $$B(d\alpha(\tilde V ), d \alpha(\tilde V ) ) =B(\tilde V ,\tilde W )$$ for all left invariant horizontal vector fields $\tilde V$ and $\tilde W $, and is thus a symplectic transformation represented in the chosen coordinates by an element of $\Sp(2n,\R)$. It follows that the group of grading preserving automorphism of the Lie algebra is  
$$\R_+ \times (\Sp(n, \R) \cup  dr(\Sp(n, \R) ) )$$
where $\R_+$ acts by $d \delta_{\lambda} $. The group of "horizontal orientation" preserving automorphisms of $\He^n$ is $\R_+ \times {\rm Sp} (2n, \R)$ and the action is this case is given by the representation
$$(\lambda, S) \to  \begin{pmatrix}
\sqrt{\lambda} S & 0 \\
0 & \lambda 
\end{pmatrix}.$$



Williamson’s Theorem  (Theorem 8.3.1, page 244, \cite{DeGosson}) gives the symplectic version of the well known diagonalisation of a positive-definite symmetric matrix by the orthogonal group. 

\begin{theorem}[Williamson] Let $M$ be a positive-definite symmetric real $2n \times 2n$ matrix. Then
	there exists $S \in \Sp(2n,\R)$ such that
		$$ S^T M S = \Lambda_M=\begin{pmatrix} \Lambda & 0 \\ 0 & \Lambda \end{pmatrix} $$ where $\Lambda= {\rm diag}(\lambda_1, \dots, \lambda_n ) $ and $\pm i \lambda_j$ are the eigenvalues of $\, -\Omega M^{-1}$.
\end{theorem}

The Williamson diagonal form gives the following factorisation 
\begin{align*}
\mathcal{L}_M u &=   S^{-1} \subh{\nabla} \cdot   \Lambda_M S^{-1} \subh{\nabla}.
\end{align*} 

Furthermore, we can dilate by $\lambda_1^{-1}$ and conclude that for any positive-definite symmetric real $2n \times 2n$ matrix, we can change coordinates by a Lie group automorphism $\alpha$ of $\He^n$ such that the frame given by $ \bar  X_i =d\alpha(\tilde X_i)$, $ \bar Y_i =d\alpha(\tilde Y_i)$ and $ \bar T = d \alpha(\tilde T)$, is left invariant and 
\begin{align}
  \mathcal{L}_M (u \circ \alpha ) &=   ( \subh{\bar \nabla} \cdot   \Lambda_M (1, \lambda_2/\lambda_1, \dots, \lambda_n/\lambda_1 )  \subh{ \bar \nabla} u) \circ \alpha \nonumber\\
   &= ( \bar X_1^2u) \circ \alpha + ( \bar Y_1^2u) \circ \alpha + \sum_{i=2}^{n} \frac{\lambda_i}{\lambda_1} (\bar  X_i^2u) \circ \alpha + \frac{\lambda_i}{\lambda_1}(\bar Y_i^2u) \circ \alpha \label{WillmForm}
\end{align} where $\subh{\bar \nabla}$ is the horizontal gradient operator relative to the frame $ \{ \bar X_i,  \bar Y_i \}$.    

\section{Sub-laplacians on $\He^2$}

In this section we consider the simplest nontrivial case of the equivalence problem for sub-Laplacians, namely $\He^2$. As a consequence of \eqref{WillmForm}, we need only to consider the equivalence problem for the one parameter family
\begin{align} \mathcal{L}_c = \tilde X_1^2+\tilde Y_1^2 +c( \tilde X_2^2+\tilde Y_2^2), \qquad 0<c. \label{Lc} \end{align} We will show that for each $c>0$, the existence of a local contact diffeomorphism $F:\He^2 \to \He^2$ with the property that $$ \mathcal{L}_c(u \circ F) = \mathcal{L}u \circ F$$ requires that the prolongation of $F$ to the second horizontal jet space $\Jet^m(\He^2, \mathbb{R})$ , restricts to a Pansu differentiable map between two stratified subgroups of $\Jet^m(\He^2, \mathbb{R})$ which implies that the Lie algebras of these two groups must be isomorphic. The nonexistence of an isomorphism between the afore mentioned Lie algebras can easily be verified with MAPLE using the package LieAlgebras. 

The main task that remains is to explicitly construct the two stratified subgroups of $\Jet^m(\He^2, \mathbb{R})$ and their Lie algebras, which can then be analysed using MAPLE. To proceed  we need a basis for $\HorDer^{\leq 2}$. By the Poincar\'e–Birkhoff–Witt theorem (\cite{MR0132805}, I.2.7), the family
\begin{align}
\{\tilde X_1^{\alpha_1}\tilde X_2^{\alpha_2}\tilde Y_1^{\beta_1}\tilde Y_2^{\beta_2}\tilde T^\gamma: \alpha_i,\beta_j, \gamma \in \N\cup \{0\}\} \label{pbwitt}
\end{align}
is a basis for the left invariant differential operators on $\He^2$. It follows that the left invariant differential operators of homogeneous order $m \leq 2$ are spanned by the set 
\begin{align} 
\{ \tilde X_1,\tilde X_2,\tilde Y_1, \tilde Y_2, \tilde X_1^2, \tilde X_1 \tilde X_2, \tilde X_1 \tilde Y_1, \tilde X_1\tilde Y_2, \tilde X_2^2, \tilde X_2 \tilde Y_1, \tilde X_2\tilde Y_2, \tilde Y_1^2, \tilde Y_1\tilde Y_2, \tilde Y_2^2, \tilde Y_1 \tilde X_1\}. \label{hset}
\end{align}
Note that we have replaced the element $\tilde T$ that is prescribed by \eqref{pbwitt} with $\tilde Y_1 \tilde X_1$.  Since  $[\tilde X_1,\tilde Y_1]=[\tilde X_2,\tilde Y_2]=-4 \tilde T$ we have     
\begin{align} \tilde T=\frac{1}{4} (\tilde Y_1 \tilde X_1 -\tilde X_1\tilde Y_1) \quad {\rm and} \quad  \tilde Y_2 \tilde X_2= \tilde X_2\tilde Y_2 + \tilde Y_1 \tilde X_1 -\tilde X_1\tilde Y_1 \label{jetnormalise} \end{align} and it follows that \eqref{hset} is a basis  for $\HorDer^{\leq 2}$. 

The two jet has the form $j^2_pu= \sum_{k=0}^2 \frac{1}{k!} A_{u,p}^k. $ Since $A_{u,p}^1(V)= \tilde V u (p)$, the set $$ \{A_{ X_1}=dx_1, A_{ X_2}=dx_2, A_{ Y_1}=dy_1,A_{ Y_2}= dy_2 \} $$ forms a basis for $\HorDer^1(\frk g; \mathbb{R})$. The second order term of $j^2_pu$ is defined on $\mathbb{V}_1^2$ by $$ A_{u,p}^2(V,W)= \tilde V  \tilde W u(p) $$ which when expressed relative to \eqref{hset} becomes
\begin{align*}A_{u,p}^2(V,W) =& \,   v_1 w_1  \tilde X_1^2u + (v_1 w_2  +  v_2 w_1)  \tilde X_1  \tilde X_2u +  (v_1 w_3 - v_4 w_2) \tilde X_1  \tilde Y_1u   \\
	&   + (v_1 w_4  + v_4 w_1)  \tilde X_1  \tilde Y_2 u + v_2 w_2  \tilde X_2^2u + (v_3 w_2 + v_2 w_3)  \tilde X_2  \tilde Y_1u    \\
	& + (v_2 w_4 + v_4 w_2) \tilde X_2  \tilde Y_2 u + v_3 w_3  \tilde Y_1^2u + (v_3 w_4 + v_4 w_3)  \tilde Y_1  \tilde Y_2u\\
	& + v_4 w_4  \tilde Y_2^2u + (v_3 w_1 + v_4 w_2)  \tilde Y_1  \tilde X_1u. 
\end{align*}  
and the basis elements of $\HorDer^2(\frk g; \mathbb{R})$ are:
\begin{align*}
 A_{ X_1^2} &=  dx_1 \otimes dx_1 \\
 A_{ X_1 X_2} &=  dx_1  \otimes dx_2 + dx_2  \otimes dx_1 \\
 A_{ X_1 Y_1} &=   dx_1  \otimes dy_1 - dy_2  \otimes dx_2  \\
 A_{ X_1 Y_2} &= dx_1  \otimes dy_2 + dy_2  \otimes dx_1\\
 A_{ X_2^2} &=   dx_2  \otimes dx_2 \\
 A_{ X_2 Y_1} &=   dy_1  \otimes dx_2 + dx_2  \otimes dy_1 \\
 A_{ X_2 Y_2} &= dx_2  \otimes dy_2 + dy_2  \otimes dx_2\\
 A_{ Y_1^2} &=  dy_1  \otimes dy_1\\
 A_{ Y_1  Y_2} &=  dy_1  \otimes dy_2 + dy_2  \otimes dy_1\\
 A_{ Y_2^2} &=   dy_2  \otimes dy_2\\
 A_{ Y_1 X_1}&=   dy_1  \otimes dx_1 + dy_2  \otimes dx_2.  
\end{align*}
  As the notation suggests, $j^2_pu= \sum_{K \in \mathcal{I} \cup \{0\} }  \tilde Ku(p) A_{K}$ where $\mathcal{I}$ is the ordered list of symbols given by $\eqref{hset}$ with tilde removed, and $\tilde Ku(p)=u(p)$ when $K=0$. 
The ordered set
\begin{align*}
E=\{ e_1, \dots e_4, A_{ X_1^2} , \dots, A_{ Y_1  X_1} , e_5, A_{ X_1}, A_{ X_2}, A_{ Y_1}, A_{ Y_2}, 1 \}
\end{align*}
forms a stratified basis for the $21$-dimensional Lie algebra $\jet^2(\frk g;\mathbb{R})$ where
\begin{align*}
\jet^2(\frk g;\mathbb{R})_1 & = \mathbb{V}_1\times \HorDer^{2}(\frk g;\mathbb{R}) = {\rm span}\, \{ e_1, \dots e_4, A_{X_1X_2}, \dots, A_{Y_1X_1}\}\\
\jet^2(\frk g;\mathbb{R})_2 & = \mathbb{V}_2\times \HorDer^{1}(\frk g;\mathbb{R}) = {\rm span}\, \{ e_5, A_{X_1}, A_{X_2}, A_{Y_1}, A_{Y_2} \}\\
\jet^2(\frk g;\mathbb{R})_3 & = \mathbb{V}_3\times \HorDer^{0}(\frk g;\mathbb{R}) = \R.
\end{align*}
Note that for $A \in \HorDer^2(\frk g; \mathbb{R})$, the choice of basis and \eqref{rcontanti} require that we define
$$[e_5, A] :=\frac{1}{4} (A(e_2,e_1) - A(e_1,e_2)).$$ 

The bracket table  for the basis $E$ is as follows (we do not include the trivial row and column generated by $E_{21}$, see \cite{J2H2JLcmw} for the MAPLE file J2H2-JLc.mw):

\noindent \resizebox{\textwidth}{!}{$ \displaystyle
	\begin {array}{ccccccccccccccccccccccc} &&&&&&&&&&&&&&&&&&&\\
	0&0& -4  E_{16} &0& -E_{17}& -E_{18}&0& -   E_{20}&0&0&0&0&0&0& -E_{19} &0& -E_{21} &0&0&0\\  
	0&0&0& -4 E_{16}&0& -E_{17} & E_{20}&0& -E_{18} & -E_{19}& -E_{20}&0&0&0& -E_{20}&0&0& -E_{21}&0&0\\  
	4 E_{16}&0&0&0&0&0& -E_{17}&0&0& -E_{18}&0& -E_{19}& -E_{20}&0&0&0&0&0& -E_{21}&0\\   
	0&4  E_{16}&0&0&0&0&0& -E_{17}&0&0& -E_{18}&0& -E_{19}& -E_{20}&0&0&0&0&0& -E_{21}\\ 
	E_{17}&0&0&0&0&0&0&0&0&0&0&0&0&0&0&0&0&0&0&0\\    
	E_{18} & E_{17}&0&0&0&0&0&0&0&0&0&0&0&0&0&0&0&0&0&0\\  
	0& -E_{20} & E_{17} &0&0&0&0&0&0&0&0&0&0&0&0& \frac{1}{4} E_{21} &0&0&0&0\\
	E_{20} &0&0 & E_{17} &0&0&0&0&0&0&0&0&0&0&0&0&0&0&0&0\\ 
	0 & E_{18} &0&0&0&0&0&0&0&0&0&0&0&0&0&0&0&0&0&0\\ 
	0 & E_{19} & E_{18 } &0&0&0&0&0&0&0&0&0&0&0&0&0&0&0&0&0\\  
	0 & E_{20} &0& E_{18} &0&0&0&0&0&0&0&0&0&0&0&0&0&0&0&0\\ 
	0&0&   E_{19} &0&0&0&0&0&0&0&0&0&0&0&0&0&0&0&0&0\\  
	0&0& E_{20} & E_{19} &0&0&0&0&0&0&0&0&0&0&0&0&0&0&0&0\\  
	0&0&0& E_{20 } &0&0&0&0&0&0&0&0&0&0&0&0&0&0&0&0\\    
	E_{19} & E_{20} &0&0&0&0&0&0&0&0&0&0&0&0&0& -\frac{1}{4} E_{21} &0&0&0&0\\  
	0&0&0&0&0&0& -\frac{1}{4} E_{21} &0&0&0&0&0&0&0& \frac{1}{4} E_{21} &0&0&0&0&0\\
	E_{21} &0&0&0&0&0&0&0&0&0&0&0&0&0&0&0&0&0&0&0\\  
	0& E_{21} &0&0&0&0&0&0&0&0&0&0&0&0&0&0&0&0&0&0\\  
	0&0& E_{21} &0&0&0&0&0&0&0&0&0&0&0&0&0&0&0&0&0\\  
	0&0&0& E_{21} &0&0&0&0&0&0&0&0&0&0&0&0&0&0&0&0\\ 
	&&&&&&&&&&&&&&&&&&&& \end {array}   $}

We denote by $\eta_k(u,p)$, the coordinates of $A_{u,p}$ relative to $E$, more specifically:
\begin{align*}
	A_{u,p}^1=\sum_{k=17}^{20} \eta_k(u,p) E_k \quad {\rm and} \quad A_{u,p}^2=\sum_{k=5}^{15} \eta_k(u,p) E_k .
\end{align*}
If $u$ is $\mathcal{L}_c$-harmonic then
\begin{align}
 \eta_5+\eta_{12}+c(\eta_9+\eta_{14})=0 \label{Lczero}
\end{align}
and we let $\Jet \mathcal{L}_c$ denote the subbundle of $\Jet^2(G, \mathbb{R}) \to G$ defined by \eqref{Lczero}. Furthermore, since $\Jet \mathcal{L}_c$ is obtained by a homogeneous linear equation of the first layer coordinates, it is a Lie subgroup of $\Jet^2(G, \mathbb{R})$. 

The bracket table for the basis of the Lie algebra of $\Jet \mathcal{L}_c$  is as follows (we do not include the trivial row and column generated by $F_{20}$, see \cite{J2H2JLcmw} for the MAPLE file JL-to-JLc.mw ):

\noindent \resizebox{\textwidth}{!}{$ \displaystyle
\begin {array}{ccccccccccccccccccc}&&&&&&&&&&&&&&&&&& \\ 0&0&-4\,F_{{15}}&0&-F_{{17}}&0&-F_{{19}}&cF_{{16}}&0&0&F_{{16}}&0&cF_{{16}}&-F_{{18}}&0&-F_{{20}}&0&0&0\\ \noalign{\medskip}0&0&0&-4\,F_{{15}}&-F_{{16}}&F_{{19}}&0&-F_{{17}}&-F_{{18}}&-F_{{19}}&0&0&0&-F_{{19}}&0&0&-F_{{20}}&0&0\\ 4\,F_{{15}}&0&0&0&0&-F_{{16}}&0&0&-F_{{17}}&0&-F_{{18}}&-F_{{19}}&0&0&0&0&0&-F_{{20}}&0\\ 0&4\,F_{{15}}&0&0&0&0&-F_{{16}}&0&0&-F_{{17}}&0&-F_{{18}}&-F_{{19}}&0&0&0&0&0&-F_{{20}}\\ F_{{17}}&F_{{16}}&0&0&0&0&0&0&0&0&0&0&0&0&0&0&0&0&0\\ 0&-F_{{19}}&F_{{16}}&0&0&0&0&0&0&0&0&0&0&0&\frac{1}{4}\,F_{{20}}&0&0&0&0\\ F_{{19}}&0&0&F_{{16}}&0&0&0&0&0&0&0&0&0&0&0&0&0&0&0\\                                                                     -cF_{{16}}&F_{{17}}&0&0&0&0&0&0&0&0&0&0&0&0&0&0&0&0&0\\  0&F_{{18}}&F_{{17}}&0&0&0&0&0&0&0&0&0&0&0&0&0&0&0&0\\                                                                                                                                                                                                                                                                                                               0&F_{{19}}&0&F_{{17}}&0&0&0&0&0&0&0&0&0&0&0&0&0&0&0\\ -F_{{16}}&0&F_{{18}}&0&0&0&0&0&0&0&0&0&0&0&0&0&0&0&0\\ 0&0&F_{{19}}&F_{{18}}&0&0&0&0&0&0&0&0&0&0&0&0&0&0&0\\ -cF_{{16}}&0&0&F_{{19}}&0&0&0&0&0&0&0&0&0&0&0&0&0&0&0\\ F_{{18}}&F_{{19}}&0&0&0&0&0&0&0&0&0&0&0&0&-\frac{1}{4}\,F_{{20}}&0&0&0&0\\ 0&0&0&0&0&-\frac{1}{4}\,F_{{20}}&0&0&0&0&0&0&0&\frac{1}{4}\,F_{{20}}&0&0&0&0&0\\ F_{{20}}&0&0&0&0&0&0&0&0&0&0&0&0&0&0&0&0&0&0\\ 0&F_{{20}}&0&0&0&0&0&0&0&0&0&0&0&0&0&0&0&0&0\\ 0&0&F_{{20}}&0&0&0&0&0&0&0&0&0&0&0&0&0&0&0&0\\ 0&0&0&F_{{20}}&0&0&0&0&0&0&0&0&0&0&0&0&0&0&0 \\ 
&&&&&&&&&&&&&&&&&& \end {array}   $}
$ $ 

Let $F$ be a local contact diffeomorphism on $\Jet^2(G, \mathbb{R})$.  If $ \mathcal{L}_c(u \circ F) = \mathcal{L}u \circ F$ then $\mathcal{L}u=0$ implies $\mathcal{L}_c(u \circ F)=0$ and so the range of $F|_ {\Jet \mathcal{L}_c}$ is $\Jet \mathcal{L} := \Jet \mathcal{L}_1$, moreover $F|_ {\Jet \mathcal{L}_c}$  maps the horizontal bundle of $\Jet \mathcal{L}_c$ to that of $\Jet \mathcal{L}$, and so it's Pansu derivative is a Lie algebra isomorphism between the Lie algebras of  $\Jet \mathcal{L}_c$ and  $\Jet \mathcal{L}$ (see section 2 and Theorem 2.9 in \cite{PansuContact} for more detail). However, by direct calculation, see \cite{J2H2JLcmw} for the MAPLE file JL-to-JLc.mw, we can show that the afore mentioned algebras are isomorphic only in the case $c=1$.

\printbibliography

\end{document}